\documentclass[12pt]{amsart}

\usepackage{latexsym,amsmath,a4wide,amssymb,amsthm,graphics}

\usepackage{comment}
\usepackage{nicematrix}
\usepackage{multirow}

\usepackage{tikz}
\usetikzlibrary{matrix}

\usepackage{arydshln}
\usepackage{setspace}

\usepackage{ytableau}
\usetikzlibrary{tikzmark}
\usepackage{nicematrix}

\usepackage{tikz-cd}

\usepackage{amsmath,lipsum}

\newcommand\undermat[2]{%
  \makebox[0pt][c]{$\smash{\underbrace{\phantom{%
\begin{matrix}#2\end{matrix}}}_{\text{$#1$}}}$}#2}

\usepackage[utf8]{inputenc}
\usepackage[OT2,T2A]{fontenc}
\usepackage{graphicx}
\usepackage{gensymb}
\usepackage{enumitem}
\usepackage{setspace}
\usepackage{tabularx}
\usepackage{upgreek}
\usepackage{hyperref}
\hypersetup{colorlinks=true, urlcolor=cyan, linkcolor=black}
\newcommand{\CC}{\mathbb{C}}

\theoremstyle{definition}
\newtheorem{remark}{Remark}[section]
\newtheorem{problem}{Problem}[section]

\newtheorem{theorem}{Theorem}[section]

\newtheorem{proposition}{Proposition}[section]
\newtheorem{example}{Example}[section]
\newtheorem{lemma}{Lemma}[section]
\newtheorem{corollary}{Corollarly}[section]

\title{On singular pencils with commuting coefficients}
\author{Vadym Koval, Patryk Pagacz}

\begin{document}
\maketitle

\begin{abstract}
We investigate the relation between the spectrum of matrix (or operator) polynomials and the Taylor spectrum of its coefficients.
We prove that the polynomial of commuting matrices is singular, i.e. its spectrum is the whole complex plane, if and only if $(0,0,\ldots,0)$ belongs to the Taylor spectrum of its coefficients. 
On the other hand we prove that this equivalence is not longer true if we consider the operators on infinite dimensional Hilbert space as coefficients of polynomial.

As a consequence we could propose a new description of (Taylor) spectrum of $k$-tuple of matrices and we could disprove the conjecture previously proposed in the literature.

Additionally, we pointed out the Kronecker forms of the pencils with commuting coefficients.

\end{abstract}


\section{Introduction}

Let $\mathcal{H}$ denote a Hilbert space and $\mathcal{B}(\mathcal{H})$ an algebra of all bounded linear transformations from $\mathcal{H}$ to $\mathcal{H}$.
If $\mathcal{H}$ is finite dimensional we can simply write $\CC^n$ instead of $\mathcal{H}$ and $\CC^{n\times n}$ instead of $\mathcal{B}(\mathcal{H})$.
By \emph{a matrix polynomial} $P(\lambda)$ we mean, a polynomial with coefficients from $\mathcal{B}(\mathcal{H})$, i.e.,
$$P(\lambda)=T_1+\lambda T_2+\lambda^2T_3+\ldots+\lambda^{k-1}T_k,\quad T_1,T_2,\ldots,T_k \in \mathcal{B}(\mathcal{H}).$$
Matrix polynomials with
degree $1$ $(k =2)$ are called \emph{matrix pencils} (or just pencils for short).

Linear pencils and, more generally, matrix polynomials find their role in applications (see e.g. \cite{NLEVP,GL,KM}). In this paper we deal with
spectra and numerical ranges connected with linear pencils and matrix polynomials.  The study of spectra and numerical ranges of matrix polynomials has a crucial role in the stability theory.
In particular, the system represented by the matrix polynomial is stable, if its spectrum is contained in the left half-plane. The
same property is shared by the numerical range. 
The classical approach to investigating matrix polynomials is linearization which reduces the eigenvalue problem for matrix polynomial to the eigenvalue problem
for linear pencils. Therefore, very often it is enough just to consider
linear pencils. However, linearizations do not transfer commutations, i.e. if we consider the matrix polynomial for commuting matrices coefficients of its linearization does not have to commute.
Therefore, if we are not able to omit the commutativity assumption, we cannot use this tool. 

In the paper \cite{KP} we focused on linear pencils with a numerical range equal to the whole complex plane. Namely, we considered the relation between linear matrix pencils with a numerical range equal to the whole complex plane and the pencils such that $(0,0)$ belongs to the joint numerical range of its coefficients. In recent paper we consider analogous questions for the spectrum of linear pencil and the Taylor spectrum of its coefficients.
A direct motivation to take up this topic was the paper \cite{DKY}, where the authors asked the question if these two classes of matrices (operators) coincide (Problems 3.17(iv),(v) in \cite{DKY}).

We start our results with 
showing that for each singular pencil $(0,0)$ belongs to the joint numerical range of its coefficients.
This theorem (Theorem \ref{dcommthm}) shows the connection between the subject of this paper and \cite{KP}. Moreover, it can be seen as an extension of the Theorem 5.4 from \cite{mehl2021matrix}.

The main part of this paper shows that the matrix polynomial is singular (its spectrum is a whole complex plane) if and only if $(0,0,\ldots,0)$ belongs to the Taylor spectrum of its coefficients (Theorem \ref{Q12_poly}). Therefore, for linear pencils, we give positive answers for Problems 3.17(iv),(v) from \cite{DKY}) in the finite dimensional case. Thanks to that we could characterize the Taylor spectrum of matrices in terms of singularity of a matrix polynomial (Corollary \ref{wn_poly}).
Another description of Taylor spectrum of invertible matrices arise in \cite{DKY} as a Conjecture 3.9. We showed that this hypothesis has a negative answer (Example \ref{hypo}). Moreover, our example is a counter-example to Theorem 3.11 from \cite{DKY}.

It has to be pointed out that the Taylor spectrum of commuting matrices coincide with some other spectra, e.g., right, left,
polynomially convex or joint algebraic spectrum, thus the choice of Taylor spectrum was important only due to infinite case. In particular, on the contrary to Theorem 1.1 from \cite{V}, the above mentioned characterization is not valid for operators on infinite dimensional Hilbert space (see Example \ref{P2}). At the end of this paper we solve Problems 3.17(i),(iv),(v) posted in \cite{DKY}, i.e. we give an example which shows that there exists a singular linear pencil (on an infinite dimension space) such that $(0,0)$ does not belong to Taylor spectrum of its coefficients (Example \ref{P2}).

Since our paper is concerned with Taylor spectrum we also discuss how the Kronecker
form looks for commuting matrices (Theorem \ref{commsing}).
From an application point of view the linear pencils with commuting coefficients are also very interesting, see \cite{KM}. Fortunately, by Lemma 2.31 from \cite{KM}, each non-singular pencil has a commuting one in its orbit.
Therefore, very often, instead of considering the pencils with commuting coefficients it is
enough to deal with non-singular pencils.
At the end of this paper we propose the conjecture (Problem \ref{Pus}) which describes all linear pencils with a commuting pencil in their orbits.

\section{Basic notions}

By \emph{a numerical range} of a matrix (or operator) polynomial $P(\lambda)$ we mean the set, $$W(P(\lambda)):=\{\lambda_0\in \CC : x^*P(\lambda_0)x=0, \textnormal{ for some nonzero } x\in \mathcal{H}\}.$$
It is easy to observe that $W(\lambda I-T)$ is equal to \emph{a numerical range of an operator} $T$, i.e., $W(T)=\{x^*Tx: x^*x=1\}$. More basic properties of the numerical ranges of matrix polynomials can be found in \cite{Li-Rodman}. 
The notion of numerical range of an operator (or a matrix) can be also extended for a $k$-tuple of operators. Namely, the set
$W(T_1,T_2,\ldots,T_k):=\{(x^*T_1x,x^*T_2x,\ldots,x^*T_kx): x^*x=1\}$ is called \emph{a joint numerical range}.

By a \emph{spectrum} of a matrix (or operator) polynomial $P(\lambda)$ we mean the set
$$\sigma(P(\lambda))=\{\lambda_0\in \CC: P(\lambda_0) \textnormal{ is not invertible}\}.$$
Thus, as for numerical ranges, $\sigma(\lambda I -T)$ is equal to $\sigma(T)$, i.e. a spectrum of an operator.
The operator polynomial (or pencil) $P(\lambda)$ such that $\sigma(P(\lambda))=\CC$ is called \emph{singular}.
For non-commuting tuples of operators there is no convenient joint spectrum. However, for commuting tuple there are several useful definitions (see \cite{MT2021}). Let us introduce a Taylor spectrum of a $k$-tuple $\textbf{T}=(T_1,T_2,\ldots,T_k)$. For this purpose we will follow \cite{Muller}.

Let $s=(s_1,\ldots,s_k)$ be a system of indeterminates. Denote by $\Lambda[s]$ the free complex algebra generated by
$s_1,\ldots,s_k$ , where the multiplication operation $\wedge$ in $\Lambda[s]$ satisfies the anticommutative relations $s_i\wedge s_j=-s_j\wedge s_i$, for $i,j=1,2,\ldots,k$. In particular, $s_i\wedge s_i=0$. For $F\subset\{1,2,\ldots,k\}$, $F=\{s_{i_1},s_{i_2},\ldots,s_{i_p}\}$, with $1\leq i_1<i_2<\ldots<i_p\leq k$ write $s_F=s_{i_1}\wedge s_{i_2}\wedge\ldots\wedge s_{i_p}$. Every element of $\Lambda[s]$ can be written uniquely in the form $\sum_{F\subset\{1,2,\ldots,k\}}\alpha_Fs_F$, with complex coefficients. Clearly, $s_{\emptyset}$ is the unit in $\Lambda[s]$.
For $p=0,1,\ldots,k$ let $\Lambda^p[s]$ be the subspace generated by the elements $s_F$, with $\# F=p$. Thus $\Lambda[s]=\bigoplus_{i=0}^k\Lambda^p[s]$.

Moreover, write $\Lambda[s,\mathcal{H}]=\mathcal{H}\otimes \Lambda[s]$ and $\Lambda^p[s,\mathcal{H}]=\mathcal{H}\otimes \Lambda^p[s]$, for a Hilbert space $\mathcal{H}$.
So
$$\Lambda[s,\mathcal{H}]=\{\sum_{F\subset\{1,2,\ldots,k\}}x_Fs_F: x_F\in \mathcal{H}\} \quad \textnormal{and } \Lambda^p[s,\mathcal{H}]=\{\sum_{F\subset\{1,2,\ldots,k\}, \#F=p}x_Fs_F: x_F\in \mathcal{H}\}.$$
Then $\Lambda[s,\mathcal{H}]$ can be considered to be also a Hilbert space, with the norm $\|\sum_{F}x_Fs_F\|=(\sum_F\|x_F\|^2)^{\frac12}$.

For $j=1,2,\ldots,k$ let $S_j\in\mathcal{B}(\Lambda[s,\mathcal{H}])$ be the operators of left multiplication by $s_j$, i.e.
$$S_j\big(\sum_Fx_Fs_F\big)=\sum_Fx_Fs_j\wedge s_F.$$

For an operator $T\in\mathcal{B}(\mathcal{H})$ we denote by the same symbol the operator
$T : \Lambda[s, \mathcal{H}] \to \Lambda[s, \mathcal{H}]$
defined by
$$T\big(\sum_Fx_Fs_F\big)=\sum_F(Tx_F)s_F.$$

Let $\textbf{T}=(T_1,T_2,\ldots,T_k)$ be an $k$-tuple of mutually commuting operators on $\mathcal{H}$.
Denote by $\delta_{\mathbf{T}}$ the operator $\delta_{\mathbf{T}}:=\sum_{i=1}^kT_iS_i\in\mathcal{B}(\Lambda[s,\mathcal{H}])$. One can check that $\delta_{\mathbf{T}}^2=0$.
Moreover, for $p=0,1,\ldots,k-1$ let $\delta^p_{\mathbf{T}}:\Lambda^p[s,\mathbf{H}]\to \Lambda^{p+1}[s,\mathbf{H}]$ be a restriction of $\delta_\mathbf{T}$ to $\Lambda^p[s,\mathbf{H}]$. Additionally, let $\delta^{-1}_{\mathbf{T}}:=0_1:\{0\}\to \Lambda^0[s,\mathcal{H}]$ and $\delta^k_{\mathbf{T}}:=0_2:\Lambda^k[s,\mathcal{H}]\to \{0\}$.

Since $\delta_{\mathbf{T}}^2=0$ the Koszul complex
$$K(\textbf{T}, \mathcal{H}) : \{0\} \xrightarrow{0_1} \Lambda^0[s,\mathcal{H}]\xrightarrow{\delta^0_{\mathbf{T}}}\Lambda^1[s,\mathcal{H}]\xrightarrow{\delta^1_{\mathbf{T}}}\ldots\xrightarrow{\delta^{k-1}_{\mathbf{T}}} \Lambda^k[s,\mathcal{H}]\xrightarrow{0_2} \{0\},$$
is well defined, i.e. $\mathcal{R}(\delta^{i-1}_{\mathbf{T}})\subset ker \delta^{i}_{\mathbf{T}}$, for $i=0,1,\ldots,k$. We say that the Koszul complex is \emph{exact} 
if $\mathcal{R}(\delta^{i-1}_{\mathbf{T}})= ker \delta^{i}_{\mathbf{T}}$, for $i=0,1,\ldots,k$.

By \emph{a Taylor spectrum} of a commuting $k$-tuple $\mathbf{T}=(T_1,T_2,\ldots,T_k)$ we mean the set
$$\sigma_T(T_1,T_2,\ldots,T_k)=\{(\lambda_1,\lambda_2\ldots,\lambda_k)\in \CC^k : K((T_1-\lambda_1,T_2-\lambda_2,\ldots,T_k-\lambda_k),\mathcal{H}) \textnormal{ is not exact} \}.$$

At the end of this section let us recall one of the main approaches to investigate the matrix pencils. Namely, a well-known Kronecker decomposition. It will be also an important tool in our considerations. 
\begin{theorem}(see \cite{G})\label{kronecker}
Let $A,B\in \CC^{n\times m}$. Then there exist invertible matrices $S\in \CC^{n\times n}$ and $T \in\CC^{m\times m}$
such that, $$S( A+\lambda B)T=diag(\mathcal{L}_{\varepsilon_1},\ldots,\mathcal{L}_{\varepsilon_k},\mathcal{L}_{\delta_1}^T,\ldots,\mathcal{L}_{\delta_r}^T,\mathcal{J}_{\gamma_1}^{\lambda_1},\ldots,\mathcal{J}_{\gamma_p}^{\lambda_p},\mathcal{N}_{\beta_1},\ldots,\mathcal{N}_{\beta_q}),$$where $\mathcal{L}_{t}\in \CC^{t\times(t+1)}$ are bidiagonal matrices of the form $$\left[\begin{array}{cccc}
0 & 1 & & \\
 & \ddots & \ddots & \\
 &  & 0 & 1 \end{array}\right]+\lambda\left[\begin{array}{cccc}
1 & 0 & & \\
 & \ddots & \ddots & \\
 &  & 1 & 0 \end{array}\right],$$
 $\mathcal{J}_{\gamma_i}^{\lambda_i}\in \CC^{\gamma_i\times\gamma_i}$ are Jordan blocks of the form
 $$\left[\begin{array}{cccc}
\lambda_i & 1 & & \\
 & \ddots &\ddots & \\
 & & \ddots & 1 \\
 &  &  & \lambda_i \end{array}\right]+\lambda\left[\begin{array}{cccc}
1 &  & & \\
 & \ddots &  & \\
 &  & \ddots & \\
 &  &  & 1 \end{array}\right],$$
 and $\mathcal{N}_{\beta_i}\in\CC^{\beta_i\times\beta_i}$ are nilpotent blocks of the form 
 $$\left[\begin{array}{cccc}
1 &  & & \\
 & \ddots &  & \\
 &  & \ddots & \\
 &  &  & 1 \end{array}\right]+\lambda\left[\begin{array}{cccc}
0 & 1 & & \\
 & \ddots &\ddots & \\
 & & \ddots & 1 \\
 &  &  & 0 \end{array}\right].$$
\end{theorem}

\section{Matrix pencils and matrix polynomials}

Before we will focus on pencils and matrix polynomials with commuting coefficients let us extend a result posted in \cite{mehl2021matrix} (Theorem 5.4).

\begin{theorem}\label{DopicoQ}
For matrices $A$, $B\in \CC^{n\times n}$, let us consider the following conditions:
\begin{itemize}
    \item[(i)] $ A+ \lambda B$ is a singular matrix pencil,
    \item[(ii)] $(0,0)\in W(A,B)$,    \item[(iii)] $W(A+\lambda B)=\CC$.
    
Then $(i)\Longrightarrow (ii) \Longrightarrow (iii)$.
\end{itemize}
\end{theorem}

\begin{proof}
$(ii) \Longrightarrow (iii)$\\
The condition $(iii)$ is  equivalent to $(0,0)\in conv W(A,B)$, see \cite{KP}.\\
$(i) \Longrightarrow (ii)$\\
Let $S$, $T\in \CC^{n\times n}$ be invertible matrices  such that
    $S( A+\lambda B)T$ has Kronecker's canonical form.
 Pencil $S( A+\lambda B)T$ is singular, therefore at least one block of the form $\mathcal{L}_{\varepsilon_i}$ and $\mathcal{L}_{\delta_j}^T$ is present. We can choose $S$ and $T$ such that $S(A+\lambda B)T$ is of the form $D=diag(\mathcal{L}_{\varepsilon_1},\mathcal{L}_{\delta_1}^T,B_1,\ldots,B_m)$, where $B_i$ are some blocks of Kronecker's canonical form. Note that from invertibility of matrices $S$ and $T$ we have 
 $(0,0)\in W(A,B)$ if and only if $(0,0)\in W(SAS^*,SBS^*)$.
 Let us denote $X=T^{-1}S^{*}$. 
 Since, $SATT^{-1}S^{*}=SAS^{*}$ and $SBTT^{-1}S^{*}=SBS^{*}$ we have to prove that the matrices of the pencil $DX$ have a common isotropic vector i.e. a vector $w$ such that $w^*SATXw=0$ and $w^*SBTXw=0$.

Without loss of generality let us assume that $D=diag(\mathcal{L}_{\varepsilon_1},\mathcal{L}_{\delta_1}^T)$ and denote $X$ as $[x_{i,j}]_{i,j=1}^{\varepsilon_1+\delta_1+1}$.
Then
$$ SATX=\begin{bmatrix}
x_{2,1} & x_{2,2} & \ldots & x_{2,\varepsilon_1+\delta_1+1}\\
    x_{3,1} & x_{3,2} & \dots & x_{3,\varepsilon_1+\delta_1+1}\\
    \vdots & \vdots & \ddots & \vdots  \\
    x_{\varepsilon_1+1,1} & x_{\varepsilon_1+1,2} & \dots & x_{\varepsilon_1+1,\varepsilon_1+\delta_1+1}\\
    0 & 0 & \ldots & 0 \\
    x_{\varepsilon_1+2,1} & x_{\varepsilon_1+2,2} & \dots & x_{\varepsilon_1+2,\varepsilon_1+\delta_1+1} \\
    \vdots & \vdots & \ddots & \vdots  \\ x_{\varepsilon_1+\delta_1+1,1} & x_{\varepsilon_1+\delta_1+1,2} & \dots & x_{\varepsilon_1+\delta_1+1,\varepsilon_1+\delta_1+1}
\end{bmatrix}$$
and 
$$ SBTX=\begin{bmatrix}
x_{1,1} & x_{1,2} & \ldots & x_{1,\varepsilon_1+\delta_1+1}\\

x_{2,1} & x_{2,2} & \ldots & x_{2,\varepsilon_1+\delta_1+1}\\
\vdots & \vdots & \ddots & \vdots  \\
    x_{\varepsilon_1,1} & x_{\varepsilon_1,2} & \dots & x_{\varepsilon_1,\varepsilon_1+\delta_1+1}\\
    x_{\varepsilon_1+2,1} & x_{\varepsilon_1+2,2} & \dots & x_{\varepsilon_1+2,\varepsilon_1+\delta_1+1}\\
    x_{\varepsilon_1+3,1} & x_{\varepsilon_1+3,2} & \dots & x_{\varepsilon_1+3,\varepsilon_1+\delta_1+1} \\
    \vdots & \vdots & \ddots & \vdots  \\ x_{\varepsilon_1+\delta_1+1,1} & x_{\varepsilon_1+\delta_1+1,2} & \dots & x_{\varepsilon_1+\delta_1+1,\varepsilon_1+\delta_1+1}\\
    0 & 0 & \ldots & 0 
\end{bmatrix}.$$

Let us consider vectors induced by the coordinates from $\varepsilon_1+1$ to $\varepsilon_1+\delta_1+1$ of the last $\delta_1$ rows of $X$ and denote them as follows
\\
$v_1=[x_{\varepsilon_1+2,\varepsilon_1+1},x_{\varepsilon_1+2,\varepsilon_1+2},\ldots,x_{\varepsilon_1+2,\varepsilon_1+\delta_1+1}],v_2=[x_{\varepsilon_1+3,\varepsilon_1+1},x_{\varepsilon_1+3,\varepsilon_1+2},\ldots,x_{\varepsilon_1+3,\varepsilon_1+\delta_1+1}],\ldots,$\\
$ v_{\delta_1}=[x_{\varepsilon_1+\delta_1+1,\varepsilon_1+1},x_{\varepsilon_1+\delta_1+1,\varepsilon_1+2},\ldots,x_{\varepsilon_1+\delta_1+1,\varepsilon_1+\delta_1+1}]$.

It is easy to see that for any $w=0^{\varepsilon_1}\times v$, where $v=(w_1,w_2,\ldots,w_{\delta_1+1})\in\CC^{\delta_1+1}$ we have $$w^*SATXw=\sum_{k=1}^{\delta_1}w_{k+1}\cdot v_k^*v \quad \textnormal{and} \quad w^*SBTXw=\sum_{k=1}^{\delta_1}w_k\cdot v_k^*v. $$
So if we choose $v\in\CC^{\delta_1+1}$ as orthogonal to $v_1,v_2,\ldots, v_{\delta_1}$ then $w$ is a common isotropic vector.

\end{proof}

In general the converse implications do not hold, even if the matrices commute.

\begin{example}{($(ii)\not\Longrightarrow (i)$)
}\label{exiiimpi}
Let $A=\begin{bmatrix}1 & 0\\ 0 & -1\end{bmatrix}
$
and $B=\begin{bmatrix}2 & 0\\ 0 & -2\end{bmatrix}
$. Then, $AB=BA, \quad A^*B=BA^*$ and $(0,0)\in W(A,B)$ ($A,B$ have a common isotropic vector), but $A+\lambda B$ is not singular. \end{example}

\begin{example}{($(iii)\not\Longrightarrow (ii)$)
}\label{exiiiimpii}
There are matrices $A,B$ such that $AB=BA$ and $W(A+\lambda  B)=\CC$, but without common isotropic vector ($(0,0)\not\in W(A,B)$). 

Recently, Lau, Li and Poon (see \cite{LauLiPoon2022}) constructed matrices $A,B\in \CC^{n\times n}$ (for $n\geq 4$) such that $W(A,B)$ is not convex and $AB=BA$.
Let $(\mu_1,\mu_2)\in conv W(A,B) \setminus W(A,B)$, then
by Theorem 2.1 from \cite{KP} the pencil $(A-\mu_1)+\lambda (B-\mu_2)$ satisfies $(iii)$, however  the matrices $A-\mu_1,B-\mu_2$ do not satisfy $(ii)$.
\end{example}

\begin{proposition}\label{dcommthm}
For matrices $A, B\in \CC^{n\times n}$ such that $AB=BA$ and $A^*B=BA^*$ the equivalence \\ $ (ii)\Longleftrightarrow (iii)$ holds.
\end{proposition}
\begin{proof}
It is known (see \cite{BR}) that a pair of doubly commuting matrices has a convex joint numerical range. Thus the implication $ (iii)\Longrightarrow (ii)$ is true due to Theorem 2.1 from \cite{KP}.
\end{proof}

For commuting matrices the condition $(i)$ in Theorem \ref{DopicoQ} can be expressed by Taylor spectrum, even more the same equivalence is true for matrix polynomials.

\begin{theorem}\label{Q12_poly}
    Let $A_1,A_2,\ldots A_k\in \CC^{n \times n}$ commute.
    Then the matrix polynomial $\sum_{i=1}^kA_i\lambda^{i-1}$ is singular if and only if $(0,0,\ldots,0)\in \sigma_T(A_1,A_2,\ldots,A_k)$.    
\end{theorem}
\begin{proof}
$(\Longleftarrow)$
        
        Let us assume that $\sum_{i=1}^kA_i\lambda^{i-1}$ is not singular. In other words, there is $\lambda_0\in \CC$ such that $\sum_{i=1}^kA_i\lambda_0^{i-1}$ is invertible. Let us note that $\big(\sum_{i=1}^kA_i\lambda_0^{i-1}\big)^{-1}$ commutes with all $A_i$.

    To show that $(0,\ldots,0)\not\in \sigma_T(A_1,\ldots,A_k)$, we will prove that the Koszul complex $K((A_1,A_2,\ldots, A_k);\mathcal{H})$ is exact. For this reason we have to prove that:
    \begin{itemize}
        \item $ker A_1\cap\ldots \cap ker A_k=ker \delta^0_{\mathbf{A}}=\mathcal{R}(\delta^{-1}_{\mathbf{A}})=\{0\}$,
        \item $\mathcal{R}(A_1)+\mathcal{R}(A_2)+\ldots+\mathcal{R}(A_k)=\mathcal{R}(\delta^{k-1}_{\mathbf{A}})=ker \delta^k_{\mathbf{A}}=\mathcal{H}$,
        \item $ker \delta^n_{\mathbf{A}} \subset \mathcal{R}(\delta_{\mathbf{A}}^{n-1})$, for $n=1,2,\ldots,k-1$.
    \end{itemize}
      The first two conditions are easy so let us fix $n$ and prove the last inclusion. 

It will be convenient to denote $\xi(i,F)=\#\{j\in F: j<i\}$, for any $F\subset\{1,2,\ldots,k\}$ and $i\in\{1,2,\ldots,k\}$. Let us assume $\sum_{F: \#F=n} x_F s_F\in ker(\delta^n)$.
Since $$\delta_{\mathbf{A}}^n(\sum_{F: \#F=n} x_F s_F)=\sum_{i=1}^k\sum_{F: \#F=n} A_ix_F s_i\wedge s_F=\sum_{F: \#F=n}\sum_{i\not\in F}(-1)^{\xi(i,F)}A_ix_F s_{F\cup \{i\}},$$
the condition $\sum_{F: \#F=n} x_F s_F\in ker(\delta^n)$ is equivalent to
$$ \sum_{i\in G}(-1)^{\xi(i,G)}A_ix_{G\setminus \{i\}}=0, $$
for any $G$ such that $\#G=n+1$.

Let us take $F$ such that $\#F=n$ and note $G_i=F\cup \{i\}$, for $i\not\in F$.
Therefore,
$$(-1)^{\xi(i,F)}A_ix_{F}+\sum_{j\in F}(-1)^{\xi(j,G_i)}A_jx_{G_i\setminus \{j\}}=\sum_{j\in G_i}(-1)^{\xi(j,G_i)}A_jx_{G_i\setminus \{j\}}=0,$$
for any $i\not\in F$.
Hence, multiplying this equality by $(-1)^{\xi(i,F)}\lambda_0^{i-1}$ and adding (for $i\not\in F$ and $i\in F$) by sides we can get,
$$\sum_{i=1}^kA_i\lambda_0^{i-1}x_F+\sum_{i\not\in F}\lambda_0^{i-1}\sum_{j\in F}(-1)^{\xi(i,F)}(-1)^{\xi(j,G_i)}A_jx_{G_i\setminus \{j\}}=\sum_{i\in F}\lambda_0^{i-1}A_ix_F.$$
Thus,
$$x_F=\sum_{j\in F}\lambda_0^{j-1}A_j \big(\sum_{l=1}^kA_l\lambda_0^{l-1}\big)^{-1}x_F-\sum_{i\not\in F}\lambda_0^{i-1}\sum_{j\in F}(-1)^{\xi(i,F)}(-1)^{\xi(j,G_i)}A_j\big(\sum_{l=1}^kA_l\lambda_0^{l-1}\big)^{-1}x_{G_i\setminus \{j\}}=$$
$$=\sum_{j\in F} 
(-1)^{\xi(j,F)} A_j \Big(
(-1)^{\xi(j,F)}\lambda_0^{j-1} Q(\lambda_0)x_F-\sum_{i\not\in F}\lambda_0^{i-1}(-1)^{\xi(j,F)}(-1)^{\xi(i,F)}(-1)^{\xi(j,G_i)}Q(\lambda_0)x_{G_i\setminus \{j\}}\Big),$$
where $Q(\lambda_0):=\big(\sum_{l=1}^kA_l\lambda_0^{l-1}\big)^{-1}$.

Let us define 
$$y_{F\setminus \{j\}}:=(-1)^{\xi(j,F)}\lambda_0^{j-1} Q(\lambda_0)x_F-\sum_{i\not\in F}\lambda_0^{i-1}(-1)^{\xi(j,F)}(-1)^{\xi(i,F)}(-1)^{\xi(j,G_i)}Q(\lambda_0)x_{G_i\setminus \{j\}},$$
for $F$ such that $\#F=n$ and $j\in F$. The sequence $y_{F\setminus \{j\}}$ is well defined. Indeed, denote $H=F\setminus\{j\}$. Since,
$$-(-1)^{\xi(j,H)}(-1)^{\xi(i,H\cup\{j\})}(-1)^{\xi(j,H\cup\{i\})}=(-1)^{\xi(i,H)}, \quad \textnormal{for } i\not\in H\cup \{j\}$$ we have that
$$y_{F\setminus \{j\}}=(-1)^{\xi(j,F)}\lambda_0^{j-1} Q(\lambda_0)x_{H\cup \{j\}}-\sum_{i\not\in H ,\ i\not=j}\lambda_0^{i-1}(-1)^{\xi(j,F)}(-1)^{\xi(i,F)}(-1)^{\xi(j,G_i)}Q(\lambda_0)x_{H\cup \{i\}}=$$
$$=(-1)^{\xi(j,H)}\lambda_0^{j-1} Q(\lambda_0)x_{H\cup \{j\}}-\sum_{i\not\in H ,\ i\not=j}\lambda_0^{i-1}(-1)^{\xi(j,H)}(-1)^{\xi(i,H\cup\{j\})}(-1)^{\xi(j,H\cup \{i\})}Q(\lambda_0)x_{H\cup \{i\}}=$$
$$=(-1)^{\xi(j,H)}\lambda_0^{j-1} Q(\lambda_0)x_{H\cup \{j\}}+\sum_{i\not\in H ,\ i\not=j}\lambda_0^{i-1}(-1)^{\xi(i,H)}Q(\lambda_0)x_{H\cup \{i\}}=$$
$$=\sum_{i\not\in H}\lambda_0^{i-1}(-1)^{\xi(i,H)}Q(\lambda_0)x_{H\cup \{i\}}.$$

Therefore,
$$\sum_{F: \#F=n}x_Fs_F= \sum_{F: \#F=n}\sum_{j\in F}(-1)^{\xi(j,F)}A_jy_{F\setminus\{j\}}s_F=\delta_{\mathbf{A}}^{n-1}\sum_{H: \#H=n-1}y_{H}s_{H}.$$
In other words, $\sum_{F: \#F=n}x_Fs_F\in \mathcal{R}(\delta^{n-1})$.

$(\Longrightarrow)$
It is well-known (see \cite{HJ}) that for commuting matrices $A_1,A_2,\ldots A_k\in\CC^{n\times n}$   there is an invertible matrix 
$S\in\CC^{n\times n}$ such that 
$$ SA_iS^{-1}=A_i^1\oplus A_i^2\oplus \ldots \oplus A_i^d, \quad \textnormal{for } i=1,2,\ldots,k,
$$
where $A_i^j$ are upper triangular with exactly one eigenvalue $\lambda_{i,j}$.
Let us assume that $\sum_{i=1}^kA_i\lambda^{i-1}$ is singular. Therefore, for any $\lambda\in \CC$ there is $j\in \{1,2,\ldots,d\}$ such that $\sum_{i=1}^k\lambda_{i,j}\lambda^{i-1}=0$. 
Hence, there is $j_0\in \{1,2,\ldots,d\}$
such that $\sum_{i=1}^k\lambda_{i,j_0}\lambda^{i-1}=0$ for infinite many $\lambda\in\CC$.
Thus $\lambda_{i,j_0}=0$, for any $i=1,2,\ldots,k$.
So $ker A_1 \cap ker A_2 \cap \ldots \cap ker A_k \not=\{0\}$ and as a consequence $(0,0,\ldots,0)\in \sigma_T(A_1,A_2,\ldots,A_k)$.
\end{proof}

Let us state the above theorem for pencils.

\begin{corollary}
    \label{Q12}
    Let $A,B\in \CC^{n \times n}$ commute.
    Then the linear pencil  $A+\lambda B$ is singular if and only if $(0,0)\in \sigma_T(A,B)$.    
\end{corollary}

A description of Taylor spectrum of matrices is a direct consequence of Theorem \ref{Q12_poly}.

\begin{corollary}\label{wn_poly}
Let $A_1, A_2, \ldots, A_k\in \CC^{n\times n}$ commute, then
$$\sigma_T(A_1,A_2,\ldots,A_k)=\{(z_1,z_2,\ldots,z_k)\in \CC^k: \textnormal{ the matrix polynomial } \sum_{i=1}^k(A_i-z_i)\lambda^{i-1} \textnormal{ is singular } \}.$$
\end{corollary}

Again let us state this result for $k=2$.

\begin{corollary}\label{wn}
Let $A, B\in \CC^{n\times n}$ commute, then
$$\sigma_T(A,B)=\{(z_1,z_2)\in \CC^2: \textnormal{ the pencil } A-z_1+\lambda(B-z_2) \textnormal{ is singular } \}.$$
\end{corollary}

We will finish this section by example which shows that Conjecture 3.9 from \cite{DKY} has negative answer and Theorem 3.11 from \cite{DKY} was not correct.

Namely, Conjecture 3.9 says that for $A,
B \in \CC^{n \times n}$ commuting invertible matrices there is
$$\sigma_T(A,B)=(\sigma(A)\times \sigma(B))\cap \{(z_1,z_2)\in \CC^2: z_1=\lambda_0 z_2 \textnormal{ for some }\lambda_0\in\sigma (A-\lambda B) \}.$$
However, it is not true even for diagonal matrices.

\begin{example}\label{hypo}
Let $A=\begin{bmatrix}
    1 & 0&0\\
    0&2&0\\
    0&0&4
\end{bmatrix}$, $B=\begin{bmatrix}
    1 & 0&0\\
    0&4&0\\
    0&0&2
\end{bmatrix}$, then $(z_1,z_2):=(2,2)\in\sigma(A)\times\sigma(B)$ and for $\lambda_0=1$ we have $z_1=\lambda_0\cdot z_2$. Moreover, $A-\lambda_0B=A-B=\begin{bmatrix}
    0&0&0\\
    0&-2&0\\
    0&0&2
\end{bmatrix}$. So $\lambda_0=1\in \sigma(A-\lambda B)$.

On the other hand, $(A-2)-(B-2)\lambda=\begin{bmatrix}
    -1+\lambda & 0&0\\
   0&-2\lambda &0\\
   0&0&2
\end{bmatrix}$ is not singular. Therefore, by Corollarly \ref{wn} we obtain $(2,2)\not\in \sigma_T(A,B)$.
\end{example}

\section{Operator pencils}


The authors of the paper 
\cite{DKY} asked the analogical question to Theorem \ref{DopicoQ}($(ii) \Longrightarrow (iii)$).
They asked if for commuting operators on Hilbert space the following implication is true:

\begin{problem}[see 3.17 (v), \cite{DKY}]\label{P1}
If $(0, 0) \in \sigma_T(T_1 , T_2)$, then is it true that $\sigma(T_1+\lambda T_2)=\CC$ ?   
\end{problem}

\begin{remark}\label{adQ12}
    It is easy to see that the first part of the proof of Theorem \ref{Q12_poly} is valid for operators. Thus Problem \ref{P1} has a positive answer, even for $k\geq 2$, i.e.
    for any $k-$tuple of commuting operators $(T_1,T_2,\ldots,T_k)$ such that $(0,0,\ldots,0)\in \sigma_T(T_1,T_2,\ldots,T_k)$, the operator polynomial $\sum_{i=1}^kT_i\lambda^{i-1}$ is singular.
\end{remark}

Therefore, the following inclusion holds.

\begin{corollary}
Let $T_1, T_2, \ldots, T_k\in \mathcal{B}(\mathcal{H})$ commute, then
$$\sigma_T(T_1,T_2,\ldots,T_k)\subset\{(z_1,z_2,\ldots,z_k)\in \CC^k: \textnormal{ the matrix polynomial } \sum_{i=1}^k(T_i-z_i)\lambda^{i-1} \textnormal{ is singular } \}.$$
\end{corollary}

In the same paper the authors asked also for the converse implication to Problem \ref{P1}:
\begin{problem}[3.17 (iv), \cite{DKY}]
Let $\sigma(T_1+\lambda T_2)=\CC$. Is it true that $(0, 0) \in\sigma_T(T_1,T_2)$?
\end{problem}

By Theorem \ref{Q12_poly} we know that the above problem has a positive answer for finite dimensional Hilbert spaces. However, the following example shows that the answer for this question, in general, is negative.

\begin{example}\label{P2}
    Let $M\in \mathcal{B}(l^2)$ be a unilateral shift, then $\sigma(M)=\mathbb{D}$. So, $M+\lambda$ is not invertible for $|\lambda|\leq 1$ and $M+\frac1\lambda$ is not invertible for $
|\lambda|\geq 1$.

    Let $T_1=I_{l^2}\oplus M \in \mathcal{B}(\mathcal{H})$ and $T_2=M\oplus I_{l^2}\in \mathcal{B}(\mathcal{H})$, where $\mathcal{H}=l^2\oplus l^2$. 

    It is easy to see that $T_1+\lambda T_2 = \lambda(\frac1\lambda I_{l^2}+ M)\oplus (M+\lambda I_{l^2})$. As a consequence $T_1+\lambda T_2$ is singular ($\sigma(T_1+\lambda T_2)=\CC$).

    Now let us see that $(0,0)\not\in \sigma_T(T_1,T_2)$. In other words, let us prove that the Koszul complex $K((T_1,T_2),\mathcal{H})$ is exact. First, let us write an equivalent form for $\delta^0_\mathbf{T}$ and $\delta^1_\mathbf{T}$:
$$\mathcal{H}\ni x\sim x s_{\emptyset}\xrightarrow{\delta^0_\mathbf{T}}T_1xs_1+T_2xs_2\sim (T_1x,T_2x)\in \mathcal{H}\oplus\mathcal{H},$$
$$\mathcal{H}\oplus\mathcal{H}\ni (x_1,x_2)\sim x_1s_1+ x_2s_2\xrightarrow{\delta^1_\mathbf{T}} (T_1x_2-T_2x_1)s_1\wedge s_2\sim T_1x_2-T_2x_1\in \mathcal{H}.$$

    Therefore, we have to show that:
    \begin{itemize}
        \item $ker T_2 \cap ker T_1=ker \delta^0_\mathbf{T}=\{0\}$,
        \item $\mathcal{R}(T_1)+ \mathcal{R}(T_2)=\mathcal{R}(\delta^1_\mathbf{T})=\mathcal{H}$,
        \item $\{(x_1,x_2): T_1x_2=T_2x_1\}=ker \delta_\mathbf{T}^1\subset \mathcal{R}(\delta^0_\mathbf{T})=\{(T_1y,T_2y): y\in \mathcal{H}\}$.
    \end{itemize}
    
    We have $ker T_2 \cap ker T_1  \subset (l^2\oplus \{0\}) \cap (\{0\}\oplus l^2) = \{0\}$. 
    Further,
    for $h=(h_1,h_2)\in l^2\oplus l^2$ we have $h=(h_1,0)+(0,h_2)=T_1(h_1,0)+T_2(0,h_2)$. So $\mathcal{R}(T_1)+ \mathcal{R}(T_2)= \mathcal{H}$.
    Finally, let us take $x_1=(x_1^1,x_1^2),x_2=(x_2^1,x_2^2)\in\mathcal{H}=l^2\oplus l^2$ such that $T_1x_2=T_2x_1$.
    Thus $$x_2^1=Mx_1^1, \quad Mx_2^2=x_1^2.$$
    So for $y=(x_1^1,x_2^2)$ we have $$(T_1y,T_2y)=(T_1(x_1^1,x_2^2),T_2(x_1^1,x_2^2))=((x_1^1,Mx_2^2),(Mx_1^1,x_2^2))=((x_1^1,x_1^2),(x_2^1,x_2^2))=(x_1,x_2).$$

Therefore, $(T_1,T_2)$ is exact and $(0,0)\not\in \sigma_T(T_1,T_2)$.
\end{example}

Moreover, the above example shows also that
the implication $(i)\Longrightarrow (ii)$ in Theorem \ref{DopicoQ} does not hold in the operator case. Indeed,
as we proved the pencil $I_{l^2}\oplus M+\lambda (M\oplus I_{l^2})$ is singular, but if $$\langle (I_{l^2}\oplus M)(h_1,h_2),(h_1,h_2) \rangle=0, \quad \langle(M\oplus I_{l^2})(h_1,h_2),(h_1,h_2)  \rangle=0$$
for some nonzero $h=(h_1,h_2)\in l^2 \oplus l^2$, then $$\|h_1\|^2=-\langle Mh_2,h_2\rangle< \|h_2\|^2=-\langle Mh_1,h_1\rangle< \|h_1\|^2.$$ 
Thus $(0,0)\not\in W(M\oplus I_{l^2},I_{l^2}\oplus M)$.

However, the second implication ($(ii)\Longrightarrow (iii)$) in Theorem \ref{DopicoQ} for operators is a direct consequence of the following extension of the Theorem 2.1 \cite{KP}.

\begin{theorem}\label{OpW}   Let $T_1,T_2\in \mathcal{B}(\mathcal{H})$, then $W(T_1+\lambda T_2)=\CC$ if and only if $(0,0)\in conv \overline{W(T_1,T_2)}$.
\end{theorem}

\begin{remark}
    The proof goes along the same line as in \cite{KP}.
    However, instead of separating the set $W(T_1,T_2)$ from zero, we have to separate the closure of $W(T_1,T_2)$.
\end{remark}

As a consequence of Theorems \ref{DopicoQ}, \ref{OpW}, Remark \ref{adQ12} and  Examples \ref{P2},\ref{exiiimpi},\ref{exiiiimpii} we can formulate the theorem analogous to Theorem \ref{DopicoQ}

\begin{theorem} 
For commuting operators $T_1,T_2\in\mathcal{B}(\mathcal{H})$ let us consider the following conditions:
\begin{itemize}
    \item[(0)] $(0,0)\in\sigma_T(T_1,T_2)$
    \item[(i)] $ \sigma(T_1+ \lambda T_2)=\CC$,
    \item[(ii)] $(0,0)\in W(T_1,T_2)$,    
    \item[(iii)] $W(T_1+\lambda T_2)=\CC$.
    
Then the following implications hold:
\[
\begin{tikzcd}
(0) \arrow{r} \arrow{d} \arrow{dr} & (i) \arrow{d} \\
(ii)  \arrow{r} & (iii).
\end{tikzcd}
\]
\end{itemize}
Moreover, any other implication does not have to hold.
\end{theorem}

\section{Kronecker forms of pencils with commuting coefficients}

The linear pencils with commuting coefficients are very useful from an application point of view, see \cite{KM}. Thus Lemma 2.31 in \cite{KM}, which says that for any non-singular linear pencil $A+\lambda B$ there exists an invertible matrix $E$ such that $EA$ and $EB$  commute, allows in many cases to relax the assumptions on the pencil $A+\lambda B$. So, often instead of considering the pencil with commuting coefficients it is enough to deal with non-singular pencils.


However, the same statement for singular pencils is not true. Before we see the possible forms of Kronecker canonical form for commuting matrices let us see the following lemma.

\begin{lemma}\label{lemmaSingComm}
    Let $A,B\in \CC^{n\times n}$ be a direct sum of canonical singular Kronecker matrices, i.e. 
    $$A+\lambda B=diag(\underbrace{\mathcal{L}_0^T,\ldots,\mathcal{L}_0^T}_{n_0},\underbrace{\mathcal{L}_{\delta_1}^T,\ldots,\mathcal{L}_{\delta_1}^T}_{n_1}\ldots,\underbrace{\mathcal{L}_{\delta_k}^T,\ldots,\mathcal{L}_{\delta_k}^T}_{n_k},$$
    \begin{equation}\label{comm_form} \underbrace{\mathcal{L}_0,\ldots,\mathcal{L}_0}_{m_0},\underbrace{\mathcal{L}_{\varepsilon_1},\ldots,\mathcal{L}_{\varepsilon_1}}_{m_1}\ldots,\underbrace{\mathcal{L}_{\varepsilon_l},\ldots,\mathcal{L}_{\varepsilon_l}}_{m_l}),
    \end{equation}where  $\mathcal{L}_{\varepsilon}$ are defined as in Theorem \ref{kronecker} with $\varepsilon_1<\varepsilon_2<\ldots<\varepsilon_l$, and $\delta_1<\delta_2<\ldots<\delta_l$.
    
    Then for any $M\in \CC^{n\times n}$ the equality $AMB=BMA$ holds if and only if $M$ has the following form
      \begin{equation}\label{Mform}
      M=\begin{array}{@{} c @{}}
\begin{bNiceArray}{c:cc|cc}[margin]
  *   & \Block{2-2}<\huge>{\mathcal{T}}      &  & 0 & 0     \\
  *  &       &  & 0 & 0    \\
\hline
* & * &*&*&* 
\\
\hdashline
 *          & \Block{2-2}<\huge>{\mathcal{A}} & & \Block{2-2}<\huge>{\mathcal{T}'} &   \\
\undermat{n_0}{*}   &    & & &                    
\end{bNiceArray}
\begin{array}{@{} r @{}}\\
      \}~\text{$m_0$}\hspace{\nulldelimiterspace} \\
      \\
    \end{array}
\end{array},\end{equation}
where
\begin{scriptsize}
\begin{equation}\label{T}\mathcal{T}=\begin{array}{@{} c @{}}
\begin{bNiceArray}[margin]{ccc|ccc|c|ccc} 
     [0]_{\delta_1\times(\delta_1+1)} & \dots & [0]_{\delta_1\times(\delta_1+1)} & [0]_{\delta_1\times(\delta_2+1)}& \dots & [0]_{\delta_1\times(\delta_2+1)}& & [0]_{\delta_1\times(\delta_k+1)} & \dots & [0]_{\delta_1\times(\delta_k+1)}\\
    
     \vdots & \ddots & \vdots & \vdots & \ddots & \vdots & \dots &\vdots & \ddots & \vdots\\    
     [0]_{\delta_1\times(\delta_1+1)} & \dots & [0]_{\delta_1\times(\delta_1+1)}&  [0]_{\delta_1\times(\delta_2+1)} & \dots & [0]_{\delta_1\times(\delta_2+1)}&  & [0]_{\delta_1\times(\delta_k+1)} & \dots & [0]_{\delta_1\times(\delta_k+1)}
    \\
    \hline
    
     \mathcal{R}_{\delta_2\times(\delta_1+1)} & \dots & \mathcal{R}_{\delta_2\times(\delta_1+1)} & [0]_{\delta_2\times(\delta_2+1)}& \dots & [0]_{\delta_2\times(\delta_2+1)}& & [0]_{\delta_2\times(\delta_k+1)} & \dots & [0]_{\delta_2\times(\delta_k+1)}\\
    
     \vdots & \ddots & \vdots & \vdots & \ddots & \vdots & \dots & \vdots & \ddots & \vdots
    \\
     \mathcal{R}_{\delta_2\times(\delta_1+1)}  & \dots & \mathcal{R}_{\delta_2\times(\delta_1+1)} & [0]_{\delta_2\times(\delta_2+1)}& \dots & [0]_{\delta_2\times(\delta_2+1)}& & [0]_{\delta_2\times(\delta_k+1)} & \dots & [0]_{\delta_2\times(\delta_k+1)}\\
    \hline
    \mathcal{R}_{\delta_3\times(\delta_1+1)} & \dots & \mathcal{R}_{\delta_3\times(\delta_1+1)} & \mathcal{R}_{\delta_3\times(\delta_2+1)}& \dots & \mathcal{R}_{\delta_3\times(\delta_2+1)}& & [0]_{\delta_3\times(\delta_k+1)} & \dots & [0]_{\delta_3\times(\delta_k+1)}\\
    
    \vdots & \ddots & \vdots & \vdots & \ddots & \vdots & \dots
    \\
    \mathcal{R}_{\delta_3\times(\delta_1+1)}  & \dots & \mathcal{R}_{\delta_3\times(\delta_1+1)} & \mathcal{R}_{\delta_3\times(\delta_2+1)}& \dots & \mathcal{R}_{\delta_3\times(\delta_2+1)}&  & [0]_{\delta_3\times(\delta_k+1)} & \dots & [0]_{\delta_3\times(\delta_k+1)}\\
    \hline
    & \vdots & & & \vdots & & \ddots & & \vdots\\
    \hline

       \mathcal{R}_{\delta_k\times(\delta_1+1)} & \dots & \mathcal{R}_{\delta_k\times(\delta_1+1)} & \mathcal{R}_{\delta_k\times(\delta_2+1)}& \dots & \mathcal{R}_{\delta_k\times(\delta_2+1)}&  & [0]_{\delta_k\times(\delta_k+1)} & \dots & [0]_{\delta_k\times(\delta_k+1)}\\
    
    \vdots & \ddots & \vdots & \vdots & \ddots & \vdots & \dots &\vdots & \ddots & \vdots
    \\
    \mathcal{R}_{\delta_k\times(\delta_1+1)}  & \dots & \mathcal{R}_{\delta_k\times(\delta_1+1)} & \mathcal{R}_{\delta_k\times(\delta_2+1)}& \dots & \mathcal{R}_{\delta_k\times(\delta_2+1)}& & [0]_{\delta_k\times(\delta_k+1)} & \dots & [0]_{\delta_k\times(\delta_k+1)}
\end{bNiceArray}
\begin{array}{@{} r @{}}
       \Bigg\}\begin{array}{@{}c@{}}\null\\\null\\\null\end{array}~\text{$n_1$}\\\null\\
       \Bigg\}\begin{array}{@{}c@{}}\null\\\null\\\null\end{array}~\text{$n_2$}\\\null\\
        \Bigg\}\begin{array}{@{}c@{}}\null\\\null\\\null\end{array}~\text{$n_3$} \\ \vdots\\\null\\\null\\\null
\end{array}
\end{array}
\end{equation}
\end{scriptsize}
with the matrices $\mathcal{R}_{\delta_i\times(\delta_j+1)}\in\CC^{\delta_i\times(\delta_j+1)}$ of the form 

\begin{equation}\label{Toeplitz}
\mathcal{R}_{\delta_i\times(\delta_j+1)}=
\begin{bmatrix}
    r_{1} & 0 & 0 & \dots & 0\\
    r_{2} & r_{1} & 0 & \dots & 0\\
    r_{3} & r_{2} & r_{1} & \dots & 0
    \\
    \vdots & \vdots & \vdots & \ddots & \vdots  \\
    r_{\delta_j+1} & r_{\delta_j} & r_{\delta_j-1} & \dots & r_{1} 
    \\
    0 & r_{\delta_j+1} & r_{\delta_j} & \dots & r_{2} 
    \\
        0 & 0 & r_{\delta_j+1} & \dots & r_{3} 
    \\

    \vdots & \vdots & \vdots & \ddots & \vdots  \\
    0 & 0 & 0 & \dots & r_{\delta_j+1} 
\end{bmatrix} 
\end{equation}
for some $r_1,r_2,\ldots,r_{\delta_j}$,\\
$\mathcal{T'}$ is a transpose of the matrix of the form \eqref{T} for $\varepsilon_1,\varepsilon_2,\ldots,\varepsilon_l $ and $m_1,m_2,\ldots,m_l$,

and
\begin{scriptsize}
\begin{equation}\label{A}\mathcal{A}=\begin{array}{@{} c @{}}
\begin{bNiceArray}[margin]{ccc|ccc|c|ccc} 
     \Lambda_{(\varepsilon_1+1)\times(\delta_1+1)} & \dots & \Lambda_{(\varepsilon_1+1)\times(\delta_1+1)} & \Lambda_{(\varepsilon_1+1)\times(\delta_2+1)}& \dots & \Lambda_{(\varepsilon_1+1)\times(\delta_2+1)}& & \Lambda_{(\varepsilon_1+1)\times(\delta_k+1)} & \dots & \Lambda_{(\varepsilon_1+1)\times(\delta_k+1)}\\
    
     \vdots & \ddots & \vdots & \vdots & \ddots & \vdots & \dots &\vdots & \ddots & \vdots\\    
     \Lambda_{(\varepsilon_1+1)\times(\delta_1+1)} & \dots & \Lambda_{(\varepsilon_1+1)\times(\delta_1+1)}&  \Lambda_{(\varepsilon_1+1)\times(\delta_2+1)} & \dots & \Lambda_{(\varepsilon_1+1)\times(\delta_2+1)}&  & \Lambda_{(\varepsilon_1+1)\times(\delta_k+1)} & \dots & \Lambda_{(\varepsilon_1+1)\times(\delta_k+1)}
    \\
    \hline
    
     \Lambda_{(\varepsilon_2+1)\times(\delta_1+1)} & \dots & \Lambda_{(\varepsilon_2+1)\times(\delta_1+1)} & \Lambda_{(\varepsilon_2+1)\times(\delta_2+1)}& \dots & \Lambda_{(\varepsilon_2+1)\times(\delta_2+1)}& & \Lambda_{(\varepsilon_2+1)\times(\delta_k+1)} & \dots & \Lambda_{(\varepsilon_2+1)\times(\delta_k+1)}\\
    
     \vdots & \ddots & \vdots & \vdots & \ddots & \vdots & \dots & \vdots & \ddots & \vdots
    \\
     \Lambda_{(\varepsilon_2+1)\times(\delta_1+1)}  & \dots & \Lambda_{(\varepsilon_2+1)\times(\delta_1+1)} & \Lambda_{(\varepsilon_2+1)\times(\delta_2+1)}& \dots & \Lambda_{(\varepsilon_2+1)\times(\delta_2+1)}& & \Lambda_{(\varepsilon_2+1)\times(\delta_k+1)} & \dots & \Lambda_{(\varepsilon_2+1)\times(\delta_k+1)}\\
    \hline
    \Lambda_{(\varepsilon_3+1)\times(\delta_1+1)} & \dots & \Lambda_{(\varepsilon_3+1)\times(\delta_1+1)} & \Lambda_{(\varepsilon_3+1)\times(\delta_2+1)}& \dots & \Lambda_{(\varepsilon_3+1)\times(\delta_2+1)}& & \Lambda_{(\varepsilon_3+1)\times(\delta_k+1)} & \dots & \Lambda_{(\varepsilon_3+1)\times(\delta_k+1)}\\
    
    \vdots & \ddots & \vdots & \vdots & \ddots & \vdots & \dots
    \\
    \Lambda_{(\varepsilon_3+1)\times(\delta_1+1)}  & \dots & \Lambda_{(\varepsilon_3+1)\times(\delta_1+1)} & \Lambda_{(\varepsilon_3+1)\times(\delta_2+1)}& \dots & \Lambda_{(\varepsilon_3+1)\times(\delta_2+1)}&  & \Lambda_{(\varepsilon_3+1)\times(\delta_k+1)} & \dots & \Lambda_{(\varepsilon_3+1)\times(\delta_k+1)}\\
    \hline
    & \vdots & & & \vdots & & \ddots & & \vdots\\
    \hline

       \Lambda_{(\varepsilon_l+1)\times(\delta_1+1)} & \dots & \Lambda_{(\varepsilon_l+1)\times(\delta_1+1)} & \Lambda_{(\varepsilon_l+1)\times(\delta_2+1)}& \dots & \Lambda_{(\varepsilon_l+1)\times(\delta_2+1)}&  & \Lambda_{(\varepsilon_l+1)\times(\delta_k+1)} & \dots & \Lambda_{(\varepsilon_l+1)\times(\delta_k+1)}\\
    
    \vdots & \ddots & \vdots & \vdots & \ddots & \vdots & \dots &\vdots & \ddots & \vdots
    \\
    \Lambda_{(\varepsilon_l+1)\times(\delta_1+1)} & \undermat{n_1}{\dots & \Lambda_{(\varepsilon_l+1)\times(\delta_1+1)} &} \Lambda_{(\varepsilon_l+1)\times(\delta_2+1)}& \undermat{n_2}{\dots & \Lambda_{(\varepsilon_l+1)\times(\delta_2+1)}&} & \Lambda_{(\varepsilon_l+1)\times(\delta_k+1)}& \undermat{n_k}{\dots & \Lambda_{(\varepsilon_l+1)\times(\delta_k+1)} }
\end{bNiceArray}
\begin{array}{@{} r @{}}
       \Bigg\}\begin{array}{@{}c@{}}\null\\\null\\\null\end{array}~\text{$m_1$}\\\null\\
       \Bigg\}\begin{array}{@{}c@{}}\null\\\null\\\null\end{array}~\text{$m_2$}\\\null\\
        \Bigg\}\begin{array}{@{}c@{}}\null\\\null\\\null\end{array}~\text{$m_3$} \\ \vdots\\\null\\\null\\\null
\end{array}
\end{array}
\end{equation}
\end{scriptsize}
\\
with the matrices $\Lambda_{(\varepsilon_i+1))\times(\delta_j+1)}\in\CC^{(\varepsilon_i+1)\times(\delta_j+1)}$ of the form 
\begin{equation}\label{Toeplitz}
\Lambda_{(\varepsilon_i+1)\times(\delta_j+1)}=
\begin{bmatrix}
    a_{1} & a_2 & a_3 & \dots & a_{\delta_j} & a_{\delta_j+1}\\
    a_{2} & a_{3} & a_4 & \dots & a_{\delta_j+1} & a_{\delta_j+2}\\
    a_{3} & a_{4} & a_{5} & \dots & a_{\delta_j+2} & a_{\delta_j+3}
    \\
    \vdots & \vdots & \vdots & \ddots & \vdots & \vdots \\
    a_{\delta_j} & a_{\delta_j+1} & a_{\delta_j+2} & \dots  & a_{\delta_j+\varepsilon_i-1} & a_{\delta_j+\varepsilon_i} \\
    a_{\delta_j+1} & a_{\delta_j+2} & a_{\delta_j+3} & \dots  & a_{\delta_j+\varepsilon_i} & a_{\delta_j+\varepsilon_i+1} 
\end{bmatrix} 
\end{equation}
for some $a_1,a_2,\ldots,a_{\delta_j+\varepsilon_i+1}$.
\end{lemma}
\begin{proof}
For simplicity's sake let us denote $(\delta_1',\delta_2',\ldots,\delta'_r):=(\underbrace{0,\ldots,0}_{n_0},\underbrace{\delta_1,\ldots,\delta_1}_{n_1}\ldots,\underbrace{\delta_k,\ldots,\delta_k}_{n_k})$ and $(\varepsilon_1',\varepsilon_2',\ldots,\varepsilon'_r):=(\underbrace{0,\ldots,0}_{m_0},\underbrace{\varepsilon_1,\ldots,\varepsilon_1}_{m_1}\ldots,\underbrace{\varepsilon_l,\ldots,\varepsilon_l}_{m_l})$.

Let $L_i^A,L_i^B\in\CC^{\varepsilon'_i\times (\varepsilon'_i+1)}$ be such that $L_i^A+\lambda L_i^B=\mathcal{L}_{\varepsilon'_i}$, for $i=1,2,\ldots ,r$,
and $L_i^A,L_i^B\in\CC^{\delta'_i\times (\delta'_i+1)}$ be such that $(L_i^{A})^T+\lambda (L_i^{B})^T=\mathcal{L}_{\delta'_i}^T$, for $i=1,2,\ldots ,r$.

For the precise forms of $L^\star_i$ see the statement of Theorem \ref{kronecker}.

Let us recall that $A$ has the form
\begin{small}
\begin{equation*}
A=\begin{bNiceArray}[margin]{ccccccccccc}

\Block{3-1}{(L_1^A)^T}&&&&\\
&&&&\\
&&&&\\
&\Block{3-1}{(L_2^A)^T}&&&\\
&&&&\\
&&&&\\
&&\Block{3-1}{(L_3^A)^T}&&\\
&&&&\\
&&&&\\
&&&\ddots &&&&\\

&&&&\Block{1-2}{L_1^A}&\\
&&&&&&\Block{1-2}{L_2^A}&\\
&&&&&&&&\Block{1-2}{L_3^A}&\\
&&&&&&&&&&\ddots \\

\CodeAfter
  \begin{tikzpicture}
  \draw (1-|1) |- (1-|2)--(4-|2)--(4-|3)--(7-|3)--(7-|4)--(10-|4)--(10-|3)--(7-|3)--(7-|2)--(4-|2)--(4-|1)--(1-|1);
  \draw[dashed] (1-|1) |- (11-|1)--(11-|5)--(1-|5)--(1-|1); 
    \draw (11-|5) |- (11-|7)--(12-|7)--(12-|9)--(13-|9)--(13-|11)--(14-|11)--(14-|9)--(13-|9)--(13-|7)--(12-|7)--(12-|5)--(11-|5);
    \draw[dashed] (11-|5) |- (11-|12)--(15-|12)--(15-|5)--(11-|5);
    
    \end{tikzpicture}
\end{bNiceArray},
\end{equation*}
\end{small}
where the first $n_0$ blocks of the form $(L^A)^T$ have degenerate dimensions $0\times 1$ and the first $m_0$ blocks of the form $L^A$ have degenerate dimensions $1\times 0$. In other words, the first $n_0$ rows are zeros and the first $m_0$ columns after the block $(L^A_r)^T$ are zeros as well.

The matrix $B$ has the same form as $A$.

Let us write $M$ in a block-matrix form, i.e.
$M=[M_{i,j}]_{i,j=1}^{2r}$, where the sizes of $M_{i,j}$ have one of the following $4$ types:
$$M_{i,j}\in\begin{cases}
\CC^{\delta'_i\times (\delta'_j+1)} \textnormal{ for }
i \leq r \textnormal{ and } j\leq r,   \\
\CC^{\delta'_i\times \varepsilon'_{j-r}} \textnormal{ for }
i \leq r \textnormal{ and } r<j\leq 2r,\\
\CC^{(\varepsilon'_{i-r}+1)\times (\delta'_j+1)} \textnormal{ for }
r<i \leq 2r \textnormal{ and } j\leq r,\\
\CC^{(\varepsilon'_{i-r}+1)\times \varepsilon'_{j-r}} \textnormal{ for }
r<i \leq 2r \textnormal{ and } r<j\leq 2r.\\
\end{cases}$$

Let us emphasize that some of the blocks $M_{i,j}$ may degenerate.

Thus $$ AMB=
\left[\begin{array}{c|c}
    [(L_i^A)^TM_{i,j}(L_j^B)^T]_{i,j=1}^r & [(L_i^A)^TM_{i,j}L_j^B]_{i=1,j=r+1}^{r,2r} \rule[-2.5ex]{0pt}{0pt}  \\
    \hline 
    [L_i^AM_{i,j}(L_j^B)^T]_{i=r+1,j=1}^{2r,r}&[L_i^AM_{i,j}L_j^B]_{i,j=r+1}^{2r} \rule[-2.5ex]{0pt}{0pt} \rule{0pt}{2.6ex}
\end{array}\right],
$$
and 
$$ BMA=
\left[\begin{array}{c|c}
    [(L_i^B)^TM_{i,j}(L_j^A)^T]_{i,j=1}^r & [(L_i^B)^TM_{i,j}L_j^A]_{i=1,j=r+1}^{r,2r} \rule[-2.5ex]{0pt}{0pt}  \\
    \hline 
    [L_i^BM_{i,j}(L_j^A)^T]_{i=r+1,j=1}^{2r,r}&[L_i^BM_{i,j}L_j^A]_{i,j=r+1}^{2r} \rule[-2.5ex]{0pt}{0pt} \rule{0pt}{2.6ex}
\end{array}\right].
$$

First, let us see the form \eqref{T}.
In other words let us prove that $M_{i,j}=0$, for $i\leq j\leq r$ and $M_{i,j}=0$, for $n_0<i,j\leq r$ such that $\delta_{i}'=\delta'_j$.

Since we have $AMB=BMA$ it is easy to see that $(L^{B}_i)^TM_{i,j}(L^{A}_j)^T=(L^{A}_i)^TM_{i,j}(L^{B}_{j})^T$, for $i,j\leq r$ and $i>n_0$. Thus for $M_{i,j}=[r_{s,t}]_{s,t=1}^{\delta'_i,\delta'_j+1}$ and $j>n_0$ (in opposite this equation gives no information) this equation has the following form:

$$\begin{bmatrix}
    0 & 0 & 0 & \dots & 0\\
    r_{1,1} & r_{1,2} & r_{1,3} & \dots & r_{1,\delta'_j}\\
    r_{2,1} & r_{2,2} & r_{2,3} & \dots & r_{2,\delta'_j}
    \\
    \vdots & \vdots & \vdots & \ddots & \vdots  \\
    r_{\delta'_i,1} & r_{\delta'_i,2} & r_{\delta'_i,3} & \dots & r_{\delta'_i,\delta'_j} 
\end{bmatrix}=\begin{bmatrix}
    r_{1,2} & r_{1,3} & r_{1,4} & \dots & r_{1,\delta'_j+1}\\
    r_{2,2} & r_{2,3} & r_{2,4} & \dots & r_{2,\delta'_j+1}
    \\
    \vdots & \vdots & \vdots & \ddots & \vdots  \\
    r_{\delta'_i,2} & r_{\delta'_i,3} & r_{\delta'_i,4} & \dots & r_{\delta'_i,\delta'_j+1} \\
    0 & 0 & 0 & \dots & 0
\end{bmatrix}.$$

Hence, if we compare the first row we get $r_{1,t}=0$ for $t\not=1$, which implies that $r_{2,t}=0$ for $t>2$. So if we compare this two matrices we will get $r_{s,t}=0$ for $t>s$.
On the other hand, if we start to compare these matrices from the bottom we would get $r_{\delta'_i-s,t}=0$ for $t<\delta'_j-s+1$.
Thus for $\delta_i'\leq \delta_j'$ we have $M_{i,j}=0$.

Moreover, for $n_0<i,j\leq r$ such that $\delta_i'>\delta'_j$, one can get that $M_{i,j}$ is a Toeplizt matrix with the following form:

\begin{equation}\label{Toeplitz}
M_{i,j}=\begin{bmatrix}
    r_{1,1} & 0 & 0 & \dots & 0\\
    r_{2,1} & r_{1,1} & 0 & \dots & 0\\
    r_{3,1} & r_{2,1} & r_{1,1} & \dots & 0
    \\
    \vdots & \vdots & \vdots & \ddots & \vdots  \\
    r_{\delta'_j+1,1} & r_{\delta'_j,1} & r_{\delta'_j-1,1} & \dots & r_{1,1} 
    \\
    0 & r_{\delta'_j+1,1} & r_{\delta'_j,1} & \dots & r_{2,1} 
    \\
        0 & 0 & r_{\delta'_j+1,1} & \dots & r_{3,1} 
    \\

    \vdots & \vdots & \vdots & \ddots & \vdots  \\
    0 & 0 & 0 & \dots & r_{\delta'_j+1,1} 
\end{bmatrix}.
\end{equation}

The same reasoning shows that $M_{i,j}=0$ for $n_0<i\leq r$ and $j>r$.

Since for $r<i,j\leq 2r$ we can consider $M_{i,j}^T$ instead of $M_{i,j}$ the properties of  $M_{i,j}$ follow from the case $i,j\leq r$.
In particular, $M_{i,j}=0$ for $i\leq r$ and for $r+m_0<i\leq 2r$ we get $M_{i,j}=0$ if $\varepsilon_i'\geq \varepsilon_j'$. Moreover, for $\varepsilon_i'< \varepsilon_j'$ the matrix $M_{i,j}$ has the form of the transpose of \eqref{Toeplitz}.

Now, let us look at $M_{i,j}$ for $r<i \leq 2r$ and $j\leq r$. We have $M_{i,j}\in \CC^{(\varepsilon'_i+1)\times (\delta'_j+1)}$
and the condition $(L_i^{A})^TM_{i,j}L_j^{B}=(L_i^{B})^TM_{i,j}L_j^{A}$ for $i>n_0$ and $j>r+m_0$ is nontrivial.
Moreover, it implies that $M_{i,j}$ without the first row and the last column is equal to $M_{i,j}$ without the first column and the last row.
So $M_{i,j}$ has the same values at antidiagonals.

\end{proof}

Now we can see the Kronecker form of the pencils with commuting coefficients.

\begin{theorem}
\label{commsing}
    Let $A,B\in\CC^{n\times n}$. If
    $A$ and $B$ commute then there are invertible matrices $S,T\in\CC^{n\times n}$ such that
    $$S(A+\lambda B)T=diag(\underbrace{\mathcal{L}_0^T,\ldots,\mathcal{L}_0^T}_{n_0},\underbrace{\mathcal{L}_{\delta_1}^T,\ldots,\mathcal{L}_{\delta_1}^T}_{n_1}\ldots,\underbrace{\mathcal{L}_{\delta_k}^T,\ldots,\mathcal{L}_{\delta_k}^T}_{n_k},$$
    \begin{equation}\label{comm_form} \underbrace{\mathcal{L}_0,\ldots,\mathcal{L}_0}_{m_0},\underbrace{\mathcal{L}_{\varepsilon_1},\ldots,\mathcal{L}_{\varepsilon_1}}_{m_1}\ldots,\underbrace{\mathcal{L}_{\varepsilon_l},\ldots,\mathcal{L}_{\varepsilon_l}}_{m_l},\mathcal{J}_{\gamma_1}^{\lambda_1},\ldots,\mathcal{J}_{\gamma_p}^{\lambda_p},\mathcal{N}_{\beta_1},\ldots,\mathcal{N}_{\beta_q}),
    \end{equation}where  $\mathcal{L}_{\varepsilon}$,$\mathcal{J}_{\gamma}^{\lambda}$, $\mathcal{N}_{\beta}$ are defined as in Theorem \ref{kronecker} with \begin{align}\label{comm_cond}  &\delta_in_i\leq n_0+n_1+\dots+n_{i-1}
    \quad \textnormal{for } i=1,\dots,k,\\
    &\varepsilon_i m_i\leq m_0+m_1+\dots+m_{i-1} \quad \textnormal{for } i=1,\dots,l.
    \nonumber
    \end{align}  
\end{theorem}

\begin{proof}
For simplicity sake let us also denote $(\delta_1',\delta_2',\ldots,\delta'_r):=(\underbrace{0,\ldots,0}_{n_0},\underbrace{\delta_1,\ldots,\delta_1}_{n_1}\ldots,\underbrace{\delta_k,\ldots,\delta_k}_{n_k})$ and $(\varepsilon_1',\varepsilon_2',\ldots,\varepsilon'_r):=(\underbrace{0,\ldots,0}_{m_0},\underbrace{\varepsilon_1,\ldots,\varepsilon_1}_{m_1}\ldots,\underbrace{\varepsilon_l,\ldots,\varepsilon_l}_{m_l})$.

Let $L_i^A,L_i^B\in\CC^{\varepsilon'_i\times (\varepsilon'_i+1)}$ be such that $L_i^A+\lambda L_i^B=\mathcal{L}_{\varepsilon'_i}$, for $i=1,2,\ldots ,r$,

$L_i^A,L_i^B\in\CC^{\delta'_i\times (\delta'_i+1)}$ be such that $(L_i^{A})^T+\lambda (L_i^{B})^T=\mathcal{L}_{\delta'_i}^T$, for $i=1,2,\ldots ,r$

$J_i^A,J_i^B\in\CC^{\gamma_i\times \gamma_i}$ be such that $J_i^A+\lambda J_i^B=\mathcal{J}_{\gamma_i}^{\lambda_i}$, for $i=1,2,\ldots ,p$,

and
$N_i^A,N_i^B\in\CC^{\beta_i\times \beta}$ be such that $N_i^A+\lambda N_i^B=\mathcal{N}_{\beta_i},$ for $i=1,2,\ldots ,q$.

The precise forms of $L^\star_i,J^\star_i,N^\star_i$ can be taken from the statement of Theorem \ref{kronecker}.

The condition
$AB=BA$ is equivalent to $SAT(T^{-1}S^{-1})SBT=SBT(T^{-1}S^{-1})SAT$. 
Let us write $(T^{-1}S^{-1})$ in block-matrix form, i.e.
$(T^{-1}S^{-1})=[R_{i,j}]_{i,j=1}^{2r+p+q}$, where the sizes of $R_{i,j}$ have one of the following $16$ types. 

$$
\left[\begin{array}{c|c|c|c}
    R_{ij}\in\CC^{\delta_i'\times (\delta_j'+1)} & R_{ij}\in\CC^{(\varepsilon_{i-r}'+1)\times (\delta_j'+1)}&R_{ij}\in\CC^{\gamma'_{i-2r}\times (\delta_j'+1)} &R_{ij}\in\CC^{\beta'_{i-p-2r}\times (\delta_j'+1)}\\ 
    \hline
    R_{ij}\in\CC^{\delta_i'\times \varepsilon'_{j-r}}&R_{ij}\in\CC^{(\varepsilon_{i-r}'+1)\times \varepsilon'_{j-r}}&R_{ij}\in\CC^{\gamma_{i-2r}'\times \varepsilon'_{j-r}}&R_{ij}\in\CC^{\beta_{i-p-2r}'\times \varepsilon'_{j-r}}\\ \hline
    R_{ij}\in\CC^{\delta_i'\times \gamma'_{j-2r}}&R_{ij}\in\CC^{(\varepsilon'_{i-r}+1)\times \gamma'_{j-2r}}&R_{ij}\in\CC^{\gamma_{i-2r}'\times \gamma'_{j-2r}}&R_{ij}\in\CC^{\beta_{i-p-2r}'\times \gamma'_{j-2r}}\\ \hline
    R_{ij}\in\CC^{\delta_i'\times \beta'_{j-p-2r}}&R_{ij}\in\CC^{(\varepsilon_{i-r}'+1)\times \beta'_{j-p-2r}}&R_{ij}\in\CC^{\gamma_{i-2r}'\times \beta'_{j-p-2r}}&R_{ij}\in\CC^{\beta_{i-p-2r}'\times \beta'_{j-p-2r}}
\end{array}\right].
$$

In other words,
$$R_{i,j}\in\begin{cases}
\CC^{\delta'_i\times (\delta'_j+1)} \textnormal{ for }
i \leq r \textnormal{ and } j\leq r,   \\
\CC^{\delta'_i\times \varepsilon'_{j-r}} \textnormal{ for }
i \leq r \textnormal{ and } r<j\leq 2r,\\
\CC^{\delta'_i\times \gamma_{j-2r}} \textnormal{ for }
i \leq r \textnormal{ and } 2r<j\leq 2r+p,\\
\CC^{\delta'_i\times \beta_{j-2r-p}} \textnormal{ for }
i \leq r \textnormal{ and } 2r+p<j,\\
\CC^{(\varepsilon'_{i-r}+1)\times (\delta'_j+1)} \textnormal{ for }
r<i \leq 2r \textnormal{ and } j\leq r,\\
\CC^{(\varepsilon'_{i-r}+1)\times \varepsilon'_{j-r}} \textnormal{ for }
r<i \leq 2r \textnormal{ and } r<j\leq 2r,\\
etc.
\end{cases}$$

Let us prove that $S^{-1}T^{-1}$ has the following form:
$$T^{-1}S^{-1}=
  \begin{array}{@{} c @{}}
\begin{bNiceArray}{c:cc|cc|c}[margin]
  *   & \Block{2-2}<\huge>{\mathcal{T}}      &  & 0 & 0 & 0    \\
  *  &       &  & 0 & 0 & 0   \\
\hline
* & * &*&*&*&* 

\\
\hdashline
 *          & \Block{2-2}<\huge>{\mathcal{A}} & & \Block{2-2}<\huge>{\mathcal{T}'} &  &* \\
*   &    & & && *                    \\
\hline
\undermat{n_0}{*} &*&* & 0 & 0 & *
\end{bNiceArray}
\begin{array}{@{} r @{}}
      \}~\text{$m_0$}\hspace{\nulldelimiterspace} \\
      \\
    \end{array}
\end{array},$$
\\

where $\mathcal{T}$ and $\mathcal{A}$ are defined as in Lemma \ref{lemmaSingComm}.
By Lemma \ref{lemmaSingComm} it is enough to show that $R_{i,j}=0$ for $i\leq r$ and $j>2r$ (the case $i> 2r$ and $r<j\leq 2r$ is just a transpose of a recent one). 

First, the equality $(L_i^A)^TR_{i,j}J_j^B=(L_i^B)^TR_{i,j}J_j^A$, for $i\leq r$ and $2r+p\geq j>2r$ has the form:
$$\left[\begin{array}{ccccc}
    0 &  0  & 0 & \cdots & 0\\
    r_{1,1} &  r_{1,2}  & r_{1,3} & \cdots & r_{1,\gamma'_{i-2r}}\\
     r_{2,1} &  r_{2,2}  & r_{2,3} & \cdots & r_{2,\gamma'_{i-2r}}\\
    \vdots &  \vdots  & \vdots & \ddots & \vdots\\
    r_{\delta'_j+1,1} &  r_{\delta'_j+1,2}  & r_{\delta'_j+1,3} & \cdots & r_{\delta'_j+1,\gamma'_{i-2r}}\\
\end{array}\right]
=$$
$$=
\left[\begin{array}{ccccc}
    \lambda_jr_{1,1} &  r_{1,1}+\lambda_jr_{1,2}  & r_{1,2}+\lambda_jr_{1,3} & \cdots & r_{1,,\gamma'_{i-2r}-1}+\lambda_jr_{1,\gamma'_{i-2r}}\\
     \lambda_jr_{2,1} & r_{2,1}+\lambda_jr_{2,2}  & r_{2,2}+\lambda_jr_{2,3} & \cdots & r_{2,\gamma'_{i-2r}-1}+\lambda_jr_{2,\gamma'_{i-2r}}\\
    \vdots &  \vdots  & \vdots & \ddots & \vdots\\
    \lambda_jr_{\delta'_j+1,1} &  r_{\delta'_j+1,1}+\lambda_jr_{\delta'_j+1,2}  & r_{\delta'_j+1,2}+\lambda_jr_{\delta'_j+1,3} & \cdots & r_{\delta'_j+1,\gamma'_{i-2r}-1}+\lambda_jr_{\delta'_j+1,\gamma'_{i-2r}}\\
    0&0&0&0&0
\end{array}\right],
$$
where $R_{i,j}=[r_{s,t}]_{s=1,t=1}^{\gamma'_{i-2r},\delta'_j+1}$.
Thus $r_{\delta'_j+1,t}=0$ for all $t=1,2,\ldots,\gamma'_{i-2r}$. And then $r_{\delta'_j,t}=0$ for all $t=1,2,\ldots,\gamma'_{i-2r}$, etc. Finally, we see that $R_{i,j}=0$.

Similarly, the equality $(L_i^A)^TR_{i,j}N_j^B=(L_i^B)^TR_{i,j}N_j^A$, for $i\leq r$ and $j>2r+p$ has the form:
$$\left[\begin{array}{ccccc}
    0 &  0  & 0 & \cdots & 0\\
   0& r_{1,1} &  r_{1,2}   & \cdots & r_{1,\beta'_{i-2r-p}-1}\\
    0& r_{2,1} &  r_{2,2}   & \cdots & r_{2,\beta'_{i-2r-p}-1}\\
    \vdots &  \vdots  & \vdots & \ddots & \vdots\\
   0& r_{\delta'_j+1,1} &  r_{\delta'_j+1,2}  & \cdots & r_{\delta'_j+1,\beta'_{i-2r-p}-1}\\
\end{array}\right]
=
\left[\begin{array}{ccccc}
    r_{1,1} &  r_{1,2}   & \cdots & r_{1,\beta'_{i-2r-p}}\\
     r_{2,1} & r_{2,2}   & \cdots & r_{2,\beta'_{i-2r-p}}\\
    \vdots &  \vdots  & \ddots & \vdots\\
    r_{\delta'_j+1,1} &  r_{\delta'_j+1,2}  & \cdots & r_{\delta'_j+1,\beta'_{i-2r-p}}\\
    0&0&0&0
\end{array}\right],
$$
where $R_{i,j}=[r_{s,t}]_{s=1,t=1}^{\beta'_{i-2r-p},\delta'_j+1}$.
Thus $r_{1,t}=0$ for all $t=1,2,\ldots,\beta'_{i-2r-p}$. And then $r_{2,t}=0$ for all $t=1,2,\ldots,\beta'_{i-2r-p}$, etc. Finally, we see that $R_{i,j}=0$.

The first $\delta_1n_1+\delta_2n_2+\dots+\delta_in_i$ rows of the matrix $T^{-1}S^{-1}$ have nonzero entries in $n_0+(\delta_1+1)n_1+(\delta_2+1)n_2+\dots+(\delta_{i-1}+1)n_{i-1}$ columns. Thus, since the rows are
linearly independent, we have
$$n_0+(\delta_1+1)n_1+(\delta_2+1)n_2+\dots+(\delta_{i-1}+1)n_{i-1}\geq \delta_1n_1+\delta_2n_2+\dots+\delta_in_i.$$
Hence, we got the first inequality in \eqref{comm_cond}. To get the second inequality we have to consider the columns which intersect the block $\mathcal{T}'$.

\end{proof}
We believe that a linear pencil $A+\lambda B$ for which there exists an invertible matrix $E$ such that $EA$ commutes with $EB$ can be characterized by the condition \eqref{comm_cond}. Thus we formulate the problem, which completes the Theorem \ref{commsing}:
\begin{problem}\label{Pus}
Let us consider the linear pencil $A+\lambda B$ with Kronecker decomposition 
   $$S(A+\lambda B)T=diag(\underbrace{\mathcal{L}_0^T,\ldots,\mathcal{L}_0^T}_{n_0},\underbrace{\mathcal{L}_{\delta_1}^T,\ldots,\mathcal{L}_{\delta_1}^T}_{n_1}\ldots,\underbrace{\mathcal{L}_{\delta_k}^T,\ldots,\mathcal{L}_{\delta_k}^T}_{n_k},$$
\begin{equation*}\underbrace{\mathcal{L}_0,\ldots,\mathcal{L}_0}_{m_0},\underbrace{\mathcal{L}_{\varepsilon_1},\ldots,\mathcal{L}_{\varepsilon_1}}_{m_1}\ldots,\underbrace{\mathcal{L}_{\varepsilon_l},\ldots,\mathcal{L}_{\varepsilon_l}}_{m_l},\mathcal{J}_{\gamma_1}^{\lambda_1},\ldots,\mathcal{J}_{\gamma_p}^{\lambda_p},\mathcal{N}_{\beta_1},\ldots,\mathcal{N}_{\beta_q}),
    \end{equation*}where $n_0<n_1<\ldots<n_k$, $m_0<m_1<\ldots<m_l$ and $\mathcal{L}_{\varepsilon}$,$\mathcal{J}_{\gamma}^{\lambda}$, $\mathcal{N}_{\beta}$ are defined as in Theorem \ref{kronecker}. 
    It is true that if the sequences $\{n_i\}_{i=1}^k$ and $\{m_j\}_{j=1}^l$ are such that 
    \begin{align}\label{Pcond}  &\delta_in_i\leq n_0+n_1+\dots+n_{i-1}
    \quad \textnormal{for } i=1,\dots,k,\\
    &\varepsilon_i m_i\leq m_0+m_1+\dots+m_{i-1} \quad \textnormal{for } i=1,\dots,l,
    \nonumber
    \end{align}
    then there is an invertible matrix $E$ such that $EA$ and $EB$ commute?
\end{problem}

By Lemma 2.31 \cite{KM} it is enough to consider Problem \ref{Pus} for the pencils with the Kronecker form consisting only of singular blocks.

Let us mention that the problem has an affirmative answer if in condition \eqref{Pcond} instead of inequalities there are equalities. First let us observe that in this case we have the following equalities
$$(\varepsilon_l+1) m_l = m_0+m_1+\dots+m_l=n_0+n_1+\dots+n_k=(\delta_k+1)n_k,$$
$$\delta_1n_1=n_0,$$
$$\delta_i n_i = n_0+n_1+\ldots+n_{i-1}=(\delta_{i-1}+1)n_{i-1}, \textnormal{ for } i=2,3,\ldots k.$$
Similarly, we have 
$$\varepsilon_1m_1=m_0, \quad \varepsilon_j m_j =(\varepsilon_{j-1}+1)m_{j-1}, \textnormal{ for } j=2,3,\ldots l.$$

Let us keep notation from Problem \ref{Pus} where $SAT=D_A, SBT=D_B$ are canonical Kronecker forms.
Now, let us define
\begin{equation*}
      E=
T^{-1}\begin{bNiceArray}{c:cc|cc}[margin]
  I_{\delta_1 n_1}   & \Block{2-2}<\huge>{\mathcal{T}}      &  & 0 & 0     \\
  0  &       &  & 0 & 0    \\
\hline
0 & 0 &0& I_{\varepsilon_1 m_1}&0 
\\
\hdashline
 0          & \Block{2-2}<\huge>{\mathcal{A}} & & \Block{2-2}<\huge>{\mathcal{T}'} &   \\
  0 &    & & &                                    
\end{bNiceArray}S^{-1},
\end{equation*}
where 
$\mathcal{T}=\begin{bmatrix}
 0&0&0&\dots &0&0\\
 I_{\delta_2 n_2}&0&0&\dots &0&0\\
 0&I_{\delta_3 n_3}&0&\dots &0&0\\
 \vdots & \vdots &\ddots & \dots & \vdots & \vdots\\
 0&0&0&\dots &I_{\delta_kn_k}&0
\end{bmatrix},
\quad 
\mathcal{T'}=\begin{bmatrix}
 0&I_{\varepsilon_2 m_2}&0&0&\dots &0\\
 0&0&I_{\varepsilon_3 m_3}&0&\dots &0\\
 \vdots & \vdots &\vdots &\ddots &\dots  & \vdots\\
 0&0&0&0&\dots &I_{\varepsilon_lm_l}\\
 0&0&0&0&\dots &0
\end{bmatrix}
$
and
$$\mathcal{A}=\begin{bmatrix}
 0&0&\dots &0&0\\
 0&0&\dots &0&0\\
 \vdots &\ddots &\dots  & \vdots& \vdots\\
 0&0&\dots &0&0\\
 0&0&\dots &0&I'_{(\delta_k+1)n_k}
\end{bmatrix}, \quad \textnormal{with } 
I'_{(\delta_k+1)n_k}=\begin{bmatrix}
    0 & \ldots & 0 & 1\\
    0 & \ldots & 1 & 0\\
    \vdots & \reflectbox{$\ddots$} & \vdots & \vdots\\
    1 & \ldots & 0 & 0
\end{bmatrix}.
$$

By the above mentioned equalities, it is easy to check that $E$ is well defined and invertible.
Moreover, since the matrix $TES$ has the form \eqref{Mform} by Lemma \ref{lemmaSingComm} we have $D_ATESD_B=D_BTESD_A$. Thus $ESD_ATESD_BT=ESD_BTESD_AT$, which is equivalent to the fact that $EA$ commutes with $EB$.

\section*{Acknowledgment}
The authors are sincerely grateful to prof. Froilan Dopico for bringing their attention to singular pencils and to prof. Micha\l \ Wojtylak for pointing them out reference \cite{KM}. 

The second author acknowledges the financial support by the JU (POB DigiWorld) grant, decision No. PSP: U1U/P06/NO/02.12

\section*{Disclosure statement}
No potential conflict of interest was reported by the authors.

\bibliographystyle{alpha}

\begin{thebibliography}{99}

\bibitem{NLEVP} T. Betcke, N. J. Higham, V. Mehrmann, C. Schr{\"o}der and
F. Tisseur. \emph{NLEVP: A collection of nonlinear eigenvalue problems}. ACM
Trans. Math. Softw., 39(2), Feb 2013.



\bibitem{DKY} S. V. Djordjević, J. Kim and J. Yoon, \emph{Spectra of
the spherical Aluthge transform, the linear pencil, and a commuting pair of operators}, Linear and
Multilinear Algebra, 70:2533-2550, 2022.

\bibitem{G} F. R. Gantmacher. Theory of Matrices, volume 1. Chelsea, New York, 1959.



\bibitem{GL} I. Gohberg, P. Lancaster and L. Rodman. \emph{Matrix polynomials}, volume 58
of Classics in Applied Mathematics. Society for Industrial and Applied Mathematics
(SIAM), Philadelphia, PA, 2009. Reprint of the 1982 original [MR0662418].


\bibitem{Harte} R.E. Harte, \emph{Invertibility, Singularity and Joseph L. Taylor}, Proc. Roy. Irish Acad. Sect. A, 81:71-79, 1981.



\bibitem{HJ} R. Horn and C. Johnson. 
Matrix Analysis. 2nd Edition,  Cambridge University Press, 2013.



\bibitem{KP} V. Koval and P. Pagacz, \emph{Matrix pencils with the numerical range equal to the whole complex plane}, Linear Algebra Appl., 657:274-286, 2023.

\bibitem{KM} P. Kunkel and V. Mehrmann, \emph{Differential-Algebraic
Equations: Analysis and Numerical Solution}, European Mathematical Society, 2006.

\bibitem{LauLiPoon2022} P.-S. Lau, C.-K. Li and Y.-T. Poon, \emph{The joint numerical range of commuting matrices}, Studia Mathematica, 267:241-259, 2022.

\bibitem{BR} V. Bolotnikov and L. Rodman, \emph{Normal forms and joint numerical ranges of doubly commuting matrices}, Linear Algebra Appl., 301:187-194, 1999.

\bibitem{mehl2021matrix} C. Mehl, V. Mehrmann and M. Wojtylak, \emph{Matrix pencils with coefficients that have positive semidefinite hermitian part}, SIAM J. Matrix Anal. Appl., 43(3):1186-1212, 2022.

\bibitem{Muller} V. M\"uller, \emph{Spectral theory of linear operators and spectral systems in {B}anach algebras},
Operator Theory: Advances and Applications, 139, Birkhäuser Verlag, Basel, 2003.

\bibitem{MT2021} V. M\"uller and Y. Tomilov, \emph{In search of convexity: diagonals and numerical ranges}, Bull. Lond. Math. Soc. 53:1016-1029, 2021. 

\bibitem{Li-Rodman} C.-K. Li and L. Rodman, \emph{Numerical range of matrix polynomials}. SIAM
J. Matrix Anal. Appl., 15(4):1256–1265, 1994.

\bibitem{V} F.H. Vasilescu, \emph{On pairs of commuting operators}, Studia Mathematica 62, 1978, pp. 203-207.

\end{thebibliography}

\address{Faculty of Mathematics and Computer Science\\
   Jagiellonian University\\
   \L ojasiewicza 6\\
   30-348 Krak\'ow\\
   Poland}
\newline
\email{vadym.koval@student.uj.edu.pl,\\patryk.pagacz@gmail.com }

\end{document}